\documentclass[12pt,leqno]{article}  

\usepackage{amsfonts}
\usepackage{amssymb}
\usepackage[all]{xy}

\newcommand{\cD}{{\cal D}}
\newcommand{\cE}{{\cal E}}

\newcommand{\cH}{{\cal H}}
\newcommand{\cI}{{\cal I}}

\newcommand{\cM}{{\cal M}}

\newcommand{\cO}{{\cal O}}
\newcommand{\cP}{{\cal P}}
\newcommand{\cQ}{{\cal Q}}
\newcommand{\cR}{{\cal R}}
\newcommand{\cV}{{\cal V}}

\newcommand{\Sym}{{\rm Sym}}
\newcommand{\Sy}{{\rm S}}
\newcommand{\rank}{{\rm rank}}

\newcommand{\coker}{{\rm coker}}
\newcommand{\ord}{{\rm ord}}

\newcommand{\Hom}{{\rm Hom}}
\newcommand{\Aut}{{\rm Aut}}
\newcommand{\id}{{\rm id}}
\newcommand{\Spec}{{\rm Spec}}

\newcommand{\Pic}{{\rm Pic}}
\newcommand{\Ind}{{\rm Ind}}

\newcommand{\Cl}{{\rm Cl}}
\newcommand{\co}{{\rm co}}
\newcommand{\lf}{{\rm lf}}

\newcommand{\ZZ}{{\mathbb Z}}

\newcommand{\NN}{{\mathbb N}}

\newcommand{\QQ}{{\mathbb Q}}

\newcommand{\wP}{{\mathfrak P}}

\renewcommand{\wp}{{\mathfrak p}}

\newcommand{\ra}{\rightarrow}
\def\rightepi{{\longrightarrow \kern-0.7em \rightarrow}}
\newcommand{\oplusm}{\mathop{\oplus}\limits}
\newcommand{\summ}{\mathop{\sum}\limits}

\newcommand{\notteilt}{{\,\not{\kern-0.075em|}\,}}
\def\antiddots{\mathinner{\mkern1mu\raise1pt\vbox{\kern7pt\hbox{.}}\mkern2mu
    \raise4pt\hbox{.}\mkern2mu\raise7pt\hbox{.}\mkern1mu}}

\begin{document}

\vspace*{15ex}

\begin{center}
{\LARGE\bf Symmetric Powers of Galois Modules on Dedekind Schemes}\\
\bigskip
by\\
\bigskip
{\sc Bernhard K\"ock}
\end{center}

\bigskip

\begin{quote}
{\footnotesize {\bf Abstract}. We prove a certain Riemann-Roch type formula
for symmetric powers of Galois modules on Dedekind schemes which, in the
number field or function field case, specializes to a formula of Burns and
Chinburg for Cassou-Nogu\`es-Taylor operations.}
\end{quote}

\bigskip

\section*{Introduction}

Let $G$ be a finite group and $E$ a number field. Let $\cO_E$ denote the ring
of integers in $E$, $Y:= \Spec(\cO_E)$, and
\[\Cl(\cO_YG) := \textrm{ker}(\rank: K_0(\cO_EG) \ra \ZZ)\]
the locally free classgroup associated with $E$ and $G$. For any $k\ge 1$,
Cassou-Nogu\`es and Taylor have constructed a certain endomorphism $\psi_k^{{\rm CNT}}$
of $\Cl(\cO_YG)$ which, via Fr\"ohlich's Hom-description of $\Cl(\cO_YG)$,
is dual to the $k$-th Adams operation on the classical ring of virtual characters
of $G$ (see \cite{CNT}). Now, let $\gcd(k, \ord(G)) =1$ and  $k'\in \NN$ an
inverse of $k$ modulo $\ord(G)$. In the paper \cite{KoCl}, we have shown that
then the endomorphism $\psi_{k'}^{{\rm CNT}}$ is a simply definable symmetric
power operation $\sigma^k$. 

Now, let $F/E$ be a finite tame Galois extension with Galois group $G$. Let
$f: X:= \Spec(\cO_F) \ra Y$ denote the corresponding $G$-morphism and $f_*$
the homomorphism
\[ f_*: K_0(G,X) \ra \Cl(\cO_YG), \quad [\cE] \mapsto 
[f_*(\cE)] - \rank(\cE) \cdot [\cO_YG],\]
from the Grothendieck group $K_0(G,X)$ of all locally free $\cO_X$-modules
with (semilinear) $G$-action to $\Cl(\cO_YG)$. Furthermore, let $\cD$ denote
the different of $F/E$ and $\psi^k$ the $k$-th Adams operation on $K_0(G,X)$.
The paper \cite{BC} by Burns and Chinburg together with the identification of
$\psi_{k'}^{{\rm CNT}}$ with $\sigma^k$ mentioned above then implies the
following Riemann-Roch type formula for all $x \in K_0(G,X)$:
\begin{equation}
\sigma^k(f_*(x)) = f_*\left(\sum_{i=0}^{k'-1} [\cD^{-ik}] \cdot \psi^k(x)\right)
\quad \textrm{in} \quad \Cl(\cO_YG)/\Ind_1^G\Cl(\cO_Y)
\end{equation}
(see Theorem 5.6 and Theorem 3.7 in \cite{KoCl}).

We now assume that $Y$ is an arbitrary Dedekind scheme and that $X$ is the
normalization of $Y$ in a finite Galois extension $F$ of the function field
$E$ of $Y$ with Galois group $G$. We again assume that the corresponding
$G$-morphism $f:X\ra Y$ is tamely ramified. Similarly to the number field case,
we define the locally free classgroup $\Cl(\cO_YG)$ (see section 2 or \cite{ABGr}),
the symmetric power operation $\sigma^k$ on $\Cl(\cO_YG)$ (see sections 1 and 2),
and the homomorphism $f_*: K_0(G,X) \ra \Cl(\cO_YG)$ (see section 3). The object
of this paper is to study the following natural question. Does the formula~(\theequation) 
still hold in this more general situation?

First of all, we mention that the paper \cite{BC} also implies that the formula 
(\theequation) holds if $Y$ is a projective smooth curve over a finite field $L$ and the
characteristic of $L$ does not divide the order of $G$ (see Theorem 3.5(b)).
In this semisimple function field case, a Hom-description of $\Cl(\cO_YG)$
again exists and the operation $\sigma^k$ is dual to the Adams operation 
$\psi^{k'}$ as in the number field case (see Theorem 2.10). In particular, 
Fr\"ohlich's techniques can be applied as in the number field case (see \cite{BC}).

In this paper, we moreover obtain the following results whose proof however
requires completely different methods since there is no Hom-description of 
$\Cl(\cO_YG)$ available in general. 

{\bf Theorem A}. The formula (\theequation) holds if one of the following assumptions is
satisfied:\\
(a) $k=1$.\\
(b) The group $G$ is Abelian and $f:X\ra Y$ is unramified.

{\bf Theorem B}. The formula (\theequation) holds after passing from 
$\Cl(\cO_YG)/\Ind_1^G\Cl(\cO_Y)$ to $\hat{K}_0(G,Y)[k^{-1}]/(\Ind_1^G K_0(Y))
\hat{K}_0(G,Y)[k^{-1}]$ via the Cartan homomorphism. \\
Here, $K_0(G,Y)$ denotes
the Grothendieck group of all locally free $\cO_Y$-modules with $G$-action
and $\hat{K}_0(G,Y)[k^{-1}]$ denotes the $I$-adic completion of $K_0(G,Y)[k^{-1}]$
where $I$ is the augmentation ideal of $K_0(G,Y)[k^{-1}]$. 

The proof of Theorem A in the case $k=1$ relies on the results of the paper \cite{Ch}
by Chase (see Proposition 3.2). Note that, despite the fact $\sigma^k = \id$ for 
$k=1$, the formula (\theequation) is non-trivial since $k'$ may be an arbitrary natural
number in the coset $1 + \ord(G)\ZZ$. If $G$ is Abelian and $f:X\ra Y$ is 
unramified, the proof of Theorem A relies on the following two facts (see
Theorem 3.5). Firstly, applying the operation $\sigma^k$ to the element
$[\cQ]-[\cP]$ in $\Cl(\cO_YG)$ is the same as pulling back the $G$-action on 
$\cP$ and $\cQ$ along the automorphism $G\ra G$, $g \mapsto g^k$ (see Theorem 2.7).
Secondly, the map $H^1(Y,G) \ra \Cl(\cO_YG)$ which maps a principal $G$-bundle
$f:X\ra Y$ to the class $[f_*(\cO_X)] - [\cO_YG]$ is a homomorphism (by Theorem
5 in the paper \cite{Wa} by Waterhouse). Theorem $B$ follows from the 
equivariant Adams-Riemann-Roch theorem (see \cite{KoGRR}) and the case $k=1$ of
Theorem A (see Theorem 3.3). Moreover, in the semisimple function field case
mentioned above, the formula (\theequation) modulo torsion can be deduced from Theorem B
if the order of $G$ is a power of a prime (see Remark 3.6).

{\bf Acknowledgments}. I would like to thank David Burns for many very helpful
discussions and for his hospitality during my stay at the King's College in
London. In particular, he has drawn my attention to the paper \cite{Ch}
which is fundamental for the first case in Theorem A and also for Theorem B.

\bigskip

\section*{\S 1 Symmetric Power Operations on $K_0$-, $K_1$-,
and Relative Grothendieck Groups}

Let $X$ be a Noetherian scheme and $G$ a finite group. 

First, we introduce the
category of locally projective modules over the group ring $\cO_XG$. Then, we
(purely algebraically) construct symmetric power operations on the Grothendieck
group $K_0(\cO_XG)$ and the Bass group $K_1^{\det}(\cO_XG)$ associated with this
category. While these constructions are more or less obvious generalizations of the
constructions in section 1 of \cite{KoCl} (for $K_0$ and $K_1$), the subsequent 
construction of symmetric power operations on relative Grothendieck groups
(in the sense of \cite{Ba}) is new. We furthermore show that these operations
are compatible with the maps in the localization sequence. Finally, we present
some cases in which the relative Grothendieck groups can be identified with 
Grothendieck groups of certain torsion modules.

By a (quasi-)coherent $\cO_XG$-module we mean a (quasi-)coherent $\cO_X$-module
$\cP$ together with an action of $G$ on $\cP$ by $\cO_X$-homomorphisms. 
Homomorphisms and exact sequences of quasi-coherent $\cO_XG$-modules are
defined in the obvious way. We call a coherent $\cO_XG$-module $\cP$ 
{\em locally projective} iff the stalk $\cP_x$ is a projective $\cO_{X,x}G$-module
for all $x\in X$. Let $K_0(\cO_XG)$ denote the Grothendieck group of all
locally projective $\cO_XG$-modules. 

{\bf Remark 1.1}. Let $X= \Spec(A)$ be affine. Then, a finitely generated module
$P$ over the group ring $AG$ is projective if and only if the corresponding
coherent $\cO_XG$-module $\cP = \tilde{P}$ is locally projective; indeed, the
exactness of the Hom-functor $\Hom_{AG}(P,-)$ is equivalent with the exactness
of the Hom-functors $\Hom_{\cO_{X,x}G}(\cP_x, -)$, $x\in X$. In particular, the
Grothendieck group $K_0(\cO_XG)$ coincides with the usual Grothendieck group
$K_0(AG)$ of all f.\ g.\ projective $AG$-modules.

We are now going to construct the above-mentioned symmetric power operations.
As in section 1 of \cite{KoCl}, it is convenient to introduce the following
categories. For any $i\ge 1$, let $\cM_i$ denote the smallest full subcategory
of the Abelian category of all coherent $\cO_XG$-modules which is closed
under extensions and kernels of $\cO_XG$-epimorphisms and which contains all
the modules of the form $\Sym_{\cO_X}^{i_1}(\cP_1) \otimes_{\cO_X} \ldots 
\otimes _{\cO_X} \Sym_{\cO_X}^{i_r}(\cP_r)$ where $\cP_1, \ldots, \cP_r$
are locally projective coherent $\cO_XG$-modules, $i_1, \ldots, i_r$ are natural
numbers with $i_1 + \ldots + i_r = i$, and $G$ acts diagonally. So, $\cM_1$ is
the category of all locally projective coherent $\cO_XG$-modules. By 
Proposition 1.1 in \cite{KoCl}, the category $\cM_i$ is contained in $\cM_1$ if
$\gcd(i, \ord(G))$ is invertible on $X$. It is easy to see that, for all $i,j 
\ge 1$, the functor 
\[\cM_i \times \cM_j \ra \cM_{i+j}, \quad (\cP, \cQ) \mapsto \cP\otimes_{\cO_X}\cQ,\]
is well-defined and bi-exact (cf.\ Lemma 1.2 in \cite{KoCl}). In particular, we
obtain products $K_0(\cM_i) \times K_0(\cM_j) \ra K_0(\cM_{i+j})$, $i,j \ge 1$,
and the set $1+\prod_{i\ge 1}K_0(\cM_i)t^i$ consisting of all power series 
$1+\sum_{i\ge 1}a_i t^i$ with $a_i \in K_0(\cM_i)$ forms an Abelian group with 
respect to multiplication of power series. As usual, one shows that
the association $[\cP] \mapsto \sum_{i\ge 0} [\Sym^i_{\cO_X}(\cP)]t^i$ can be 
extended to a well-defined homomorphism
\[\sigma: K_0(\cO_XG) \ra 1+\prod_{i\ge 1}K_0(\cM_i) t^i\]
(see \S 1 of Chapter V in \cite{FL} and Lemma 1.3 in \cite{KoCl}). The
$i$-th component of this homomorphism is denoted by $\sigma^i$. We have for
all $x,y \in K_0(\cO_XG)$: 
\begin{eqnarray*}
\sigma^i(x-y) &=& \sum_{{a\ge 0, b_1, \ldots, b_u \ge 1} \atop {a+b_1+\ldots+b_u =i}}
(-1)^u \sigma^a(x) \sigma^{b_1}(y) \cdots \sigma^{b_u}(y) \\
&=& \sum_{{a,b_1, \ldots, b_u \ge 1} \atop {a+b_1+\ldots+b_u =i}}
(-1)^u (\sigma^a(x)-\sigma^a(y))\sigma^{b_1}(y) \cdots \sigma^{b_u}(y)
\end{eqnarray*}
in $K_0(\cM_i)$ (cf.\ section 2 in \cite{Gr2}). If $\gcd(i,\ord(G))$ is invertible
on $X$, let $\sigma^i$ also denote the composition
\[K_0(\cO_XG) \,\, \stackrel{\sigma^i}{\longrightarrow}\,\,
K_0(\cM_i) \,\, \stackrel{{\rm can}}{\longrightarrow}\,\, K_0(\cO_XG).\]
The map $\sigma^i$ is called {\em $i$-th symmetric power operation}.

Now, let $K_0(\ZZ,\cM_i)$ denote the Grothendieck group of all pairs 
$(\cP,\alpha)$ where $\cP$ is an object of $\cM_i$ and $\alpha$ is an 
$\cO_XG$-automorphism of $\cP$. We put $K_0(\ZZ,\cO_XG):=K_0(\ZZ,\cM_1)$. As 
above, the association $((\cP,\alpha), (\cQ,\beta)) \mapsto (\cP \otimes_{\cO_X}
\cQ, \alpha \otimes_{\cO_X} \beta)$ induces a multiplication map
\[K_0(\ZZ, \cM_i) \times K_0(\ZZ, \cM_j) \ra K_0(\ZZ,\cM_{i+j})\]
(for all $i, j \ge 1$) and the association $(\cP,\alpha) \mapsto 
\sum_{i\ge 0}(\Sym^i_{\cO_X}(\cP), \Sym^i_{\cO_X}(\alpha))t^i$ induces 
a homomorphism
\[\sigma: K_0(\ZZ,\cO_XG) \ra 1+\prod_{i\ge 1}K_0(\ZZ, \cM_i)t^i.\]
By restricting, we obtain symmetric power operations 
\[\sigma^i: \tilde{K}_0(\ZZ,\cO_XG)=\tilde{K}_0(\ZZ, \cM_1) \ra 
\tilde{K}_0(\ZZ,\cM_i),\quad i\ge1,\]
between the reduced Grothendieck groups
\[\tilde{K}_0(\ZZ,\cM_i) := \textrm{ker}(K_0(\ZZ,\cM_i) \ra K_0(\cM_i), \;
[\cP,\alpha] \mapsto [\cP]), \quad i\ge1.\]
We denote the factor group of $K_0(\ZZ,\cM_i)$ modulo the subgroup generated
by the relations of the form $[\cP,\alpha\beta]-[\cP,\alpha]-[\cP,\beta]$ by
$K_1^{\det}(\cM_i)$. If $X=\Spec(A)$ is affine, the group $K_1^{\det}(\cO_XG) =
K_1^{\det}(\cM_1)$ coincides with the usual Bass-Whitehead group $K_1(AG)$
of the group ring $AG$ (by Remark 1.1). In the sequel, 
we consider $K_1^{\det}(\cM_i)$ as the
factor group of $\tilde{K}_0(\ZZ,\cM_i)$ modulo the subgroup $I_i$ generated by
the relations of the form $[\cP,\alpha\beta] - [\cP,\alpha] - [\cP,\beta] + [\cP,\id]$.
Since 
\begin{eqnarray*}
\lefteqn{([\cP,\alpha\beta] - [\cP,\alpha] - [\cP,\beta] + [\cP,\id]) \cdot
[\cQ,\gamma]}\\
&=& \left([\cP\otimes\cQ, \alpha\beta\otimes\gamma] - [\cP\otimes\cQ, \alpha\otimes\gamma]
- [\cP \otimes \cQ, \beta \otimes \id] + [\cP\otimes\cQ, \id \otimes \id]\right) \\
&& - \left([\cP\otimes\cQ, \beta\otimes\gamma] - [\cP\otimes\cQ,\id\otimes\gamma]
-[\cP\otimes\cQ, \beta\otimes\id] + [\cP\otimes\cQ, \id \otimes \id] \right),
\end{eqnarray*}
the group $I_iK_0(\ZZ,\cM_j)$ is contained in $I_{i+j}$ and we obtain a 
multiplication map
\[K_1^{\det}(\cM_i) \times K_1^{\det}(\cM_j) = \tilde{K}_0(\ZZ,\cM_i)/I_i \times
\tilde{K}_0(\ZZ,\cM_j)/I_j \ra K_1^{\det}(\cM_{i+j})\]
(for all $i,j \ge 1$) which is obviously trivial, i.e., the product of any
two power series $\sum_{i\ge 0}x_it^i$, $\sum_{i\ge 0} y_it^i$ in 
$1 + \prod_{i\ge 1}K_1^{\det}(\cM_i)t^i$ is $1+ \sum_{i\ge 1}(x_i + y_i)t_i$.

{\bf Lemma 1.2}. The homomorphism $\sigma: \tilde{K}_0(\ZZ, \cO_XG) \ra
1+\prod_{i\ge 1}\tilde{K}_0(\ZZ,\cM_i)t^i$ induces a homomorphism
$\sigma: K_1^{\det}(\cO_XG) \ra 1 + \prod_{i\ge 1} K_1^{\det}(\cM_i)t^i$.
Each component $\sigma^i: K_1^{\det}(\cO_XG) \ra K_1^{\det}(\cM_i)$ of $\sigma$ 
is a homomorphism. 

{\bf Proof}. Let $\cP \in \cM_1$ and $\alpha, \beta \in \Aut_{\cO_XG}(\cP)$. 
We write $\Sy$ for $\Sym$. Then, for all $a \ge 1$, the element
\begin{eqnarray*}
\lefteqn{[\Sy^a(\cP\oplus\cP), \Sy^a(\alpha\beta \oplus \id)] - 
[\Sy^a(\cP\oplus\cP), \Sy^a(\alpha \oplus \beta)]   }\\
&=&\sum_{c=0}^a \Big([\Sy^c(\cP)\otimes \Sy^{a-c}(\cP), 
\Sy^c(\alpha\beta) \otimes \Sy^{a-c}(\id)]
- [\Sy^c(\cP) \otimes \Sy^{a-c}(\cP), \Sy^c(\alpha) \otimes \id]\\
&&-[\Sy^c(\cP)\otimes \Sy^{a-c}(\cP), \Sy^c(\beta) \otimes \id] +
[\Sy^c(\cP)\otimes\Sy^{a-c}(\cP), \id\otimes \id]\Big)\\
&& - \sum_{c=0}^a \Big([\Sy^c(\cP) \otimes \Sy^{a-c}(\cP), \Sy^c(\alpha) \otimes
\Sy^{a-c}(\beta)]
- [\Sy^c(\cP) \otimes \Sy^{a-c}(\cP), \Sy^c(\alpha) \otimes \id]\\
&&-[\Sy^c(\cP) \otimes \Sy^{a-c}(\cP), \id \otimes \Sy^{a-c}(\beta)]
+[\Sy^c(\cP) \otimes \Sy^{a-c}(\cP), \id\otimes \id]\Big)
\end{eqnarray*}
is contained in $I_a$. Since 
\[\sigma^i(x-y) = \sum_{{a,b_1, \ldots, b_u \ge 1} \atop {a+b_1+\ldots+b_u =i}}
(-1)^u(\sigma^a(x) - \sigma^a(y)) \sigma^{b_1}(y) \cdots \sigma^{b_u}(y)\]
(for all $x,y \in K_0(\ZZ, \cM_1)$), this implies that the element
\[\sigma^i\left([\cP,\alpha\beta]-[\cP,\alpha]-[\cP,\beta]+[\cP,\id]\right)
=\sigma^i\left([\cP\oplus \cP,\alpha \beta \oplus \id] - [\cP\oplus \cP,\alpha \oplus \beta]\right)\]
is contained in $I_i$, as was to be shown. For all $x,y \in K_1^{\det}(\cM_1)$,
we have
\[\sigma(x+y) = \sigma(x) \cdot \sigma(y) = 1+\sum_{i\ge 1}
(\sigma^i(x) + \sigma^i(y))t^i \quad \textrm{in} \quad 1+\prod_{i\ge 1} 
K_1^{\det}(\cM_i)t^i;\]
thus, $\sigma^i$ is a homomorphism for all $i\ge 1$. 

Now, let $j: U \ra X$ be a morphism between Noetherian schemes. Similarly to
\S 5 of Chapter VII in \cite{Ba}, let $K_0(\co(j_i^*))$ denote the Grothendieck
group of all triples $(\cP, \alpha, \cQ)$ where $\cP$ and $\cQ$ are objects 
in $\cM_i$ and $\alpha: j^*(\cP) \ra j^*(\cQ)$ is an $\cO_UG$-isomorphism. 
As above, the association 
\[((\cP,\alpha, \cQ),(\cP',\alpha',\cQ')) \mapsto
(\cP \otimes_{\cO_X}\cP', \alpha\otimes_{\cO_U}\alpha', \cQ\otimes_{\cO_X}\cQ')\]
induces, for all $i,i'\ge 1$, a multiplication map
\[K_0(\co(j_i^*)) \times K_0(\co(j_{i'}^*)) \ra K_0(\co(j_{i+i'}^*))\]
and the association $(\cP, \alpha, \cQ) \mapsto \sum_{i\ge 0}
(\Sym^i_{\cO_X}(\cP), \Sym^i_{\cO_U}(\alpha), \Sym^i_{\cO_X}(\cQ))t^i$
induces a homomorphism
\[\sigma: K_0(\co(j_1^*)) \ra 1+ \prod_{i\ge 1} K_0(\co(j_i^*))t^i.\]
By restricting, we obtain symmetric power operations
\[\sigma^i: \tilde{K}_0(\co(j_1^*)) \ra \tilde{K}_0(\co(j_i^*)), \quad i\ge 1,\]
between the reduced Grothendieck groups
\[\tilde{K}_0(\co(j_i^*)):= \textrm{ker}(K_0(\co(j_i^*)) \ra K_0(\cM_i), \;
[\cP,\alpha,\cQ] \mapsto [\cP]).\]
Let $K_0(j_i^*)$ denote the factor group of $K_0(\co(j_i^*))$ modulo the 
subgroup generated by the relations of the form $[\cP, \beta \alpha, \cR]
-[\cP, \alpha, \cQ] - [\cQ, \beta, \cR]$ (see also Proposition (5.1) on p.\ 370
in \cite{Ba}). In the sequel, we consider $K_0(j_i^*)$ as the factor group
of $\tilde{K}_0(\co(j_i^*))$ modulo the subgroup $I_i$ generated by the
elements of the form $[\cP,\beta \alpha, \cR] - [\cP, \alpha, \cQ] 
- [\cQ, \beta, \cR] + [\cQ, \id, \cQ]$. As above, one easily sees that
$I_i K_0(\co(j_{i'}^*))$ is contained in $I_{i+i'}$
and we obtain a multiplication map
\[K_0(j_i^*) \times K_0(j_{i'}^*) \ra K_0(j_{i+i'}^*)\]
for all $i,i' \ge 1$ which however (in contrast to $K_1^{\det}$) seems not 
to be trivial in general. 

{\bf Lemma 1.3}. The homomorphism $\sigma: \tilde{K}_0(\co(j_1^*)) \ra
1+ \prod _{i\ge 1} \tilde{K}_0(\co(j_i^*))t^i$ induces a homomorphism
$\sigma: K_0(j_1^*) \ra 1+ \prod_{i\ge 1} K_0(j_i^*) t^i$. 

{\bf Proof}. Similarly to Lemma 1.2.

The association $[\cP, \alpha, \cQ] \mapsto [\cQ] - [\cP]$ obviously defines 
a homomorphism 
\[\nu_i: K_0(j_i^*) \ra K_0(\cM_i)\]
for all $i\ge 1$. 

{\bf Lemma 1.4}. The multiplication maps are compatible with the homomorphisms
$\nu_i$, $i\ge 1$. The same holds for the symmetric power operations $\sigma^i$,
$i\ge 1$; i.e., the following diagram commutes for all $i\ge 1$:
\[\xymatrix{
K_0(j_1^*) \ar[r]^{\nu_1} \ar[d]^{\sigma^i}& K_0(\cM_1) \ar[d]^{\sigma^i}\\
K_0(j_i^*) \ar[r]^{\nu_i} & K_0(\cM_i). }\]

{\bf Proof}. We only prove the assertion for $\sigma^i$. Let $\cP, \cQ, \cR \in \cM_1$
and $\alpha: j^*(\cP) \,\, \tilde{\ra}\,\, j^*(\cQ)$, $\beta: j^*(\cQ) \,\,
\tilde{\ra} \,\, j^*(\cR)$ $\cO_UG$-isomorphisms. We again write $\Sy$ for $\Sym$.
Then we have in $K_0(\cM_i)$:
\begin{eqnarray*}
\lefteqn{\nu_i\sigma^i(\cP,\alpha,\cQ) = \nu_i\sigma^i((\cP,\alpha,\cQ)-(\cP,\id,\cP))}\\
&=& \nu_i\Bigg(\sum_{{a \ge 0, b_1, \ldots, b_u \ge 1} \atop {a+b_1+\ldots+b_u = i}}
(-1)^u \Big(\Sy^a(\cP) \otimes \Sy^{b_1}(\cP)\otimes \ldots \otimes \Sy^{b_u}(\cP),
\Sy^a(\alpha) \otimes \id \otimes \ldots \otimes \id, \\
&& \hspace*{43ex}\Sy^a(\cQ) \otimes \Sy^{b_1}(\cP) \otimes 
\ldots \otimes \Sy^{b_u}(\cP)\Big)\Bigg)\\
&=& \sum_{{a,b_1, \ldots, b_u \ge 1}\atop {a+b_1+\ldots+ b_u=i}}(-1)^u
\left([\Sy^a(\cQ)] - [\Sy^a(\cP)]\right) \cdot [\Sy^{b_1}(\cP)\otimes \ldots 
\otimes \Sy^{b_u}(\cP)]\\
&=&\sigma^i([\cQ]-[\cP]) = \sigma^i\nu_1(\cP,\alpha,\cQ).
\end{eqnarray*}

We now assume that $U = \Spec(F)$ is affine. Then, by Proposition (2.1) on p.\ 393 in
\cite{Ba}, the association $(\oplusm^m FG,\alpha) \mapsto 
(\oplusm^m\cO_XG,\alpha,\oplusm^m \cO_XG)$ induces a {\em connecting homomorphism}
\[\partial: K_1(FG)\ra K_0(j_1^*)\]
with $\nu_1 \circ \partial=0$.

{\bf Lemma 1.5}. Let $\gcd(i, \ord(G))$ be invertible on $X$. Then we have:
\[\sigma^i \circ \partial = \partial \circ \sigma^i \quad \textrm{in} \quad
\Hom(K_1(FG),K_0(j_1^*)).\]
The multiplication maps are compatible with $\partial$ (in the obvious sense), too.
In particular, the multiplication on $\textrm{Image}(\partial)$ is trivial and
the operation $\sigma^i$ is a homomorphism on $\textrm{Image}(\partial)$. 

{\bf Proof}. Easy.

{\bf Proposition 1.6}. The following sequence is exact:
\[K_1(FG) \,\, \stackrel{\partial}{\ra} \,\, K_0(j_1^*) \,\,
\stackrel{\nu_1}{\ra} \,\, K_0(\cO_XG)\,\, \stackrel{j^*}{\ra}\,\, K_0(FG).\]

{\bf Proof}. Apply Theorem (2.2)(b) on p.\ 396 in \cite{Ba}.

Now, let $\cH$ denote the category of all coherent $\cO_XG$-modules $\cV$ which
allow a resolution by locally projective coherent $\cO_XG$-modules of length
$\le 1$ and for which $j^*(\cV)=0$ holds. Furthermore, let $K_0T(\cO_XG)$ denote the
Grothendieck group of $\cH$. By mapping the class $[\cV]$ of a coherent
$\cO_XG$-module $\cV$ with the resolution $0 \ra \cP \,\, \stackrel{\alpha}{\ra}\,\,
\cQ \ra \cV \ra 0$ and with $j^*(\cV) =0$ to the element $(\cP, j^*(\alpha), \cQ)$
in $K_0(j_1^*)$, we obviously obtain a homomorphism 
\[\psi: K_0T(\cO_XG) \ra K_0(j_1^*).\] 

{\bf Proposition 1.7}. The homomorphism $\psi$ is bijective in the following cases:\\
(a) $X = \Spec(A)$ is affine, $F$ is the localization $A_S$ of $A$ by a multiplicative
set $S$ of non-zero-divisors in $A$, and $j: U=\Spec(F) \ra X=\Spec(A)$ is the 
canonical morphism.\\
(b) The morphism $j:U=\Spec(F) \ra X$ is an open immersion and the ideal $\cI$
of the complement $Y:=X\backslash U$ is locally generated by a non-zero-divisor.\\
(c) $X$ is a Dedekind scheme (i.e., Noetherian, regular, irreducible, and 
$\dim(X) =1$), $F$ is the function field of $X$ and $j:U=\Spec(F) \ra X$ is the
canonical morphism.

{\bf Proof}. The assertion (a) follows from (the proof of) Theorem (5.8) on p.\ 429
in \cite{Ba}. In the case~(b), we construct an 
inverse map as follows: Let $(\cP,\alpha,\cQ)$
be a generator of $K_0(j_1^*)$. Then, the image of the composition
\[\tilde{\alpha}: \cP \,\,\stackrel{{\rm can}}{\longrightarrow}\,\, j_*j^*(\cP) \,\,
\stackrel{j_*(\alpha)}{\longrightarrow} \,\, j_*j^*(\cQ) = 
\cup_{n\ge 0} \cI^{-n}\cQ \]
(see Lemma 2 on p.\ 231 in \cite{Gr1} for the last equality) is contained in 
$\cI^{-n}\cQ$ for some $n\ge 0$. We put 
\[\phi(\cP,\alpha,\cQ):= [\coker(\cP\,\,\stackrel{\tilde{\alpha}}{\hookrightarrow}\,\,
\cI^{-n}\cQ)]-[\coker(\cQ\,\,\stackrel{{\rm can}}{\hookrightarrow}\,\, 
\cI^{-n}\cQ)] \in K_0T(\cO_XG).\]
As in loc.\ cit., one easily checks that the association $(\cP,\alpha,\cQ) \mapsto
\phi(\cP,\alpha,\cQ)$ induces a well-defined map $\phi: K_0(j_1^*)\ra K_0T(\cO_XG)$
which is an inverse of $\psi$. In the case~(c), we construct an inverse map as
follows. Let $(\cP,\alpha,\cQ)$ be a generator of $K_0(j_1^*)$. The isomorphism
$\alpha: j^*(\cP) \,\, \tilde{\ra}\,\, j^*(\cQ)$ can be extended to an isomorphism
$\cP|_U \,\, \tilde{\ra}\,\, \cQ|_U$ where $U$ is an open subset of $X$. The ideal
$\cI$ of the complement $Y:= X\backslash U$ is then locally generated by a
non-zero-divisor. We now define $\phi(\cP,\alpha,\cQ)$ as in the case (b). As in
loc.\ cit., one again easily checks that the association $(\cP, \alpha, \cQ)
\mapsto \phi(\cP,\alpha,\cQ)$ induces a well-defined map $\phi: K_0(j_1^*) \ra
K_0T(\cO_XG)$ which is an inverse of $\psi$. 

{\bf Remark 1.8}. We assume that one of the conditions (a), (b), (c) of 
Proposition 1.7 holds.\\
(a) The $K$-theory space of the exact category $\cH$ is homotopy 
equivalent to the homotopy fibre of the canonical continuous map from the
$K$-theory space of $\cM_1$ to the $K$-theory space of the exact category 
consisting of all f.\ g.\ projective $FG$-modules (see \cite{Gr1} and \cite{ABGr}).
Hence, we have a long exact (localization) sequence
\[\ldots \ra K_1(FG) \ra K_0T(\cO_XG) \ra K_0(\cO_XG) \ra K_0(FG).\]
The end of this sequence can be identified with the exact sequence in Proposition 1.6
by virtue of Proposition 1.7.\\
(b) If $\gcd(i,\ord(G))$ is invertible on $X$, we obtain a symmetric power operation
$\sigma^i:K_0T(\cO_XG) \ra K_0T(\cO_XG)$ by virtue of the isomorphism $\psi$. It maps
the class $[\cV]$ of a coherent $\cO_XG$-module $\cV$ in $\cH$ with the resolution
$0\ra \cP \,\, \stackrel{\alpha}{\ra} \,\, \cQ \ra \cV \ra 0$ to the element
\begin{eqnarray*}
\lefteqn{
\summ_{{a, b_1, \ldots, b_u\ge 1} \atop {a+b_1+\ldots + b_u=i}} (-1)^u
\Big[\coker\Big(\Sym^a(\cP) \otimes 
\Sym^{b_1}(\cP) \otimes \ldots \otimes \Sym^{b_u}(\cP)}\\ 
&&\xymatrix{ {}\ar[rrr]^{\Sym^a(\alpha)\otimes \id \otimes \ldots \otimes \id} &&&{}}
\Sym^a(\cQ) \otimes \Sym^{b_1}(\cP) \otimes \ldots \otimes \Sym^{b_u}(\cP)\Big)\Big]. 
\end{eqnarray*}
Alternatively, the operation $\sigma^i$ on $K_0T(\cO_XG)$ 
can also be constructed as follows.
Let $\cE$ denote the exact category of all short exact sequences 
$0\ra \cP \ra \cQ \ra \cV \ra 0$ with $\cP, \cQ \in \cM_1$ and $\cV \in \cH$. 
Then, we have a canonical isomorphism
\[K_0T(\cO_XG) = K_0(\cH) \cong \textrm{ker}(K_0(\cE) \ra K_0(\cO_XG), \;
[0\ra \cP \ra \cQ \ra \cV \ra 0] \mapsto [\cP]).\]
The association 
\[[0\ra \cP\,\, \stackrel{\alpha}{\ra}\,\, \cQ \ra \cV\ra 0]
\mapsto [0 \ra \Sym^i(\cP) \,\, \stackrel{\Sym^i(\alpha)}{\longrightarrow}\,\,
\Sym^i(\cQ) \ra \coker(\Sym^i(\alpha))\ra 0]\] 
induces an operation $\sigma^i$ 
on $K_0(\cE)$ as usual. It is then easy to check that its restriction
to $K_0T(\cO_XG)$ coincides with the operation $\sigma^i$ constructed above. 
Moreover, the latter construction can be extended to all higher $K$-groups 
$K_q(\cH)$, $q\ge 0$, by using the methods of \cite{Gr2}. On the other hand,
we have a symmetric power operation $\sigma^i$ on the $K$-theory space of $\cM_1$
and on the $K$-theory space of the category consisting of all f.\ g.\ projective
modules (see section 1 in \cite{KoCl}), hence also on the homotopy fibre 
mentioned in (a) and finally on $K_q(\cH)$, $q\ge 0$. It seems to be plausible
that these two constructions of $\sigma^i$ on $K_q(\cH)$, $q \ge 0$, coincide.
I hope to say more on this in a future paper.

\bigskip

\section*{\S 2 Symmetric Power Operations on Locally Free Classgroups of
Dedekind Schemes}

Let $X$ be a Dedekind scheme (i.e., Noetherian, regular, irreducible and
$\dim(X) \le 1$) with function field $F$, and let $G$ be a finite group.

First, we recall the definition of the locally free classgroup 
$\Cl(\cO_XG)$ (see \cite{ABGr}
or \cite{BC}). Using the tools developed in section 1 and Hattori's theorem,
we then show that the locally free classgroup coincides with the analogously
defined locally projective classgroup and that the operations $\sigma^i$, $i\ge 1$,
constructed in section 1 are homomorphisms on $\Cl(\cO_XG)$. Furthermore, we
prove the following concrete interpretations of the operations $\sigma^i$, $i\ge 1$,
on $\Cl(\cO_XG)$. Firstly, if $G$ is Abelian and $\gcd(i,\ord(G))=1$, then 
pulling back the action of $G$ on locally free $\cO_XG$-modules along
the automorphism $G\ra G, \; g\mapsto g^i$, induces the operation $\sigma^i$ on
$\Cl(\cO_XG)$. Secondly, if $X$ is a smooth curve over an 
(algebraically closed or) finite field $L$ such that the characteristic of $L$ does not
divide the order of $G$, then the identification of the locally free with
the locally projective classgroup allows us a simple module theoretic 
description of the isomorphism between $\Cl(\cO_XG)$ and $\Hom_{{\rm Galois}}
(K_0(\bar{L}G),\Cl(\bar{X}))$ (developed in \cite{ABGr}), and the operation
$\sigma^i$ on $\Cl(\cO_XG)$ is dual to the adjoint Adams operation $\hat{\psi}^i$
on $K_0(\bar{L}G)$ with respect to this isomorphism. The proof of the latter
result presented here can also be applied in the number field case and then
simplifies the proof of Theorem 3.7 in \cite{KoCl}.

A coherent $\cO_XG$-module $\cP$ is called {\em locally free over $\cO_XG$} iff
the stalk $\cP_x$ is a free $\cO_{X,x}G$-module for all $x \in X$. By 
Proposition (30.17) on p.\ 627 in \cite{CR}, this is equivalent to the condition
that $\cP_x \otimes_{\cO_{X,x}}\hat{\cO}_{X,x}$ is a free 
$\hat{\cO}_{X,x}G$-module for all closed points $x\in X$. (Here, $\hat{\cO}_{X,x}$
denotes the ${\mathfrak m}_x$-adic completion of $\cO_{X,x}$ and ${\mathfrak m}_x$
the maximal ideal in $\cO_{X,x}$.) Let $K_0^\lf(\cO_XG)$ denote the Grothendieck
group of all coherent $\cO_XG$-modules which are locally free over $\cO_XG$.

{\bf Remark 2.1}. Let $X=\Spec(A)$ be affine. Then we also write $K_0^\lf(AG)$
for $K_0^\lf(\cO_XG)$. This is the Grothendieck group considered for instance in
\cite{F1}. If $A$ is a local Dedekind domain, then the rank (over $AG$) induces
an isomorphism $K_0^\lf(AG) \,\, \tilde{\ra}\,\, \ZZ$. If ${\rm char}(A) =0$ and
and no prime divisor of $\ord(G)$ is a unit in $A$, then any f.\ g.\
projective $AG$-module is already locally free by Swan's theorem (see Theorem
(32.11) on p.\ 676 in \cite{CR}). The same holds if $p={\rm char}(A) > 0$ and
$G$ is a $p$-group since then the group rings $\cO_{X,x}G$, $x\in X$, are
local rings. We will prove in Proposition 2.4 that the locally free
classgroup defined below always coincides with the analogously defined locally projective
classgroup. 

{\bf Definition 2.2}. The group
\[\Cl(\cO_XG) := \textrm{ker}(K_0^\lf(\cO_XG) \,\, \stackrel{{\rm can}}{\longrightarrow}
\,\, K_0^\lf(FG) \cong \ZZ)\]
is called the {\em locally free classgroup associated with $X$ and $G$}.

Let $K_0T(\cO_XG)$ (resp., $K_0^\lf T(\cO_XG)$) denote the Grothendieck group
of all coherent $\cO_XG$-modules which are $\cO_X$-torsion modules and which
allow a resolution of length $\le 1$ by locally projective (resp., locally free)
$\cO_XG$-modules. The notation $K_0T(\cO_XG)$ obviously agrees with the
notation introduced in section 1 (if $j:U=\Spec(F) \ra X$ is the canonical
morphism). 

{\bf Lemma 2.3}. The canonical homomorphisms
\[K_0T(\cO_XG) \ra \oplusm_{x \in X \;{\rm closed}} K_0T(\cO_{X,x}G) \quad 
\textrm{and} \quad  K_0^\lf T(\cO_XG) \ra \oplusm_{x\in X \;{\rm closed}}
K_0^\lf T(\cO_{X,x}G)\]
are bijective.

{\bf Proof}. Let $x$ be a closed point of $X$ and $V$ a f.\ g.\ $\cO_{X,x}G$-module
which is $\cO_{X,x}$-torsion and which allows an $\cO_{X,x}G$-projective 
(resp., $\cO_{X,x}G$-free) resolution $0\ra P \ra Q \,\, \stackrel{\varepsilon}{\ra} 
V \ra 0$. Let $i: \Spec(\cO_{X,x}) \hookrightarrow X$ denote the inclusion. It suffices to
show that $i_*(V)$ has a (global) locally projective (resp., locally free) 
resolution of length $\le 1$. If $P$ and $Q$ are $\cO_{X,x}G$-free, i.e., if
they are isomorphic to $\oplusm^m \cO_{X,x}G$ for some $m \ge 0$, then the 
composition $\tilde{\varepsilon}: \oplusm^m \cO_XG \,\, 
\stackrel{{\rm can}}{\longrightarrow} \,\, i_*(\oplusm^m(\cO_{X,x}G)) \,\,
\stackrel{i_*(\varepsilon)}{\longrightarrow} \,\, i_*(V)$ is surjective and 
$\textrm{ker}(\varepsilon)$ is a locally free $\cO_XG$-module, i.e., $i_*(V)$ has
a locally free resolution of length $1$. If $P$ and $Q$ are only projective over
$\cO_{X,x}G$, we choose a (non-equivariant) surjective homomorphism $\cE \ra
i_*(V)$ with a locally free $\cO_X$-module $\cE$. Then, the induced 
homomorphism $\tilde{\varepsilon}: \cO_XG \otimes_{\cO_X} \cE \ra
i_*(V)$ is an equivariant surjection and the coherent $\cO_XG$-module
$\textrm{ker}(\varepsilon)$ is locally projective by Schanuel's Lemma, i.e.,
$i_*(V)$ has a locally projective resolution of length $1$. 

{\bf Proposition 2.4}. The canonical homomorphism $K_0^\lf(\cO_XG) \ra
K_0(\cO_XG)$ induces an isomorphism
\[\Cl(\cO_XG) \,\, \tilde{\ra}\,\, \textrm{ker}(K_0(\cO_XG) \,\,
\stackrel{{\rm can}}{\longrightarrow} \,\, K_0(FG)).\]

{\bf Proof}. We have a natural commutative diagram of groups
\[\xymatrix{
K_1(FG) \ar[r] \ar@{=}[d] & K_0^\lf T(\cO_XG) \ar[r] \ar[d] & 
K_0^\lf(\cO_XG) \ar[r] \ar[d] & K_0^\lf(FG) \ar[d] \\
K_1(FG) \ar[r] & K_0T(\cO_XG) \ar[r] & K_0(\cO_XG) \ar[r] & K_0(FG); }\]
here, the lower row is the exact localization sequence constructed in 
Proposition 1.6 and Proposition 1.7; the maps in the upper row are defined
as in the lower row; one can prove as in section 1 or as in Theorem 1(ii) on 
p.\ 3 in \cite{F2} that also the upper sequence is exact. Thus, it suffices
to prove that the map $K_0^\lf T(\cO_XG) \ra K_0T(\cO_XG)$ is bijective. By
Lemma 2.3, it furthermore suffices to prove that the map $K_0^\lf T(\cO_{X,x}G) 
\ra K_0T(\cO_{X,x}G)$ is bijective for all closed points $x\in X$. We
have a natural commutative diagram of groups 
\[\xymatrix{
K_1(\cO_{X,x}G) \ar[r] \ar@{=}[d] & K_1(FG) \ar[r] \ar@{=}[d] &
K_0^\lf T(\cO_{X,x}G) \ar[r] \ar[d] & 0 \\
K_1(\cO_{X,x}G) \ar[r] &  K_1(FG) \ar[r]& K_0T(\cO_{X,x}G) \ar[r] & 
K_0(\cO_{X,x}G) \ar[r] & K_0(FG) }\]
with exact rows (e.g., see Theorem 1(ii) on p.\ 3 in \cite{F2}). Furthermore, the 
map $K_0(\cO_{X,x}G) \ra K_0(FG)$ is injective by Hattori's Theorem (see
Theorem (32.1) on p.\ 671 in \cite{CR}). This proves Proposition 2.4.

Let $K_0(G,X)$ denote the Grothendieck group of all coherent $\cO_XG$-modules
which are locally free as $\cO_X$-modules.

{\bf Corollary 2.5}. If $\ord(G)$ is invertible on $X$, the Cartan
homomorphism $K_0^\lf(\cO_XG) \ra K_0(G,\cO_X)$ induces an isomorphism
\[\Cl(\cO_XG) \,\, \tilde{\ra} \,\, \textrm{ker}\left(K_0(G,X) \,\,
\stackrel{{\rm can}}{\longrightarrow} \,\, K_0(G,F) \cong K_0(FG)\right).\]

{\bf Proof}. This immediately follows from Proposition 2.4 and the fact that
a f.\ g.\ $\cO_{X,x}G$-module is projective over $\cO_{X,x}G$ if and only if
it is projective over $\cO_{X,x}$. 

Now, we fix $i\in \NN$ such that $\gcd(i,\ord(G))$ is invertible on $X$. 
By section 1, we have a symmetric power operation $\sigma^i:
K_0(\cO_XG) \ra K_0(\cO_XG)$. By restricting, we obtain an operation
$\sigma^i$ on $\textrm{ker}(K_0(\cO_XG) \ra K_0(FG)) \cong \Cl(\cO_XG)$. 
In the same way, we obtain a multiplication map on $\Cl(\cO_XG)$. 

{\bf Proposition 2.6}. The multiplication on $\Cl(\cO_XG)$ is trivial and
the operation $\sigma^i$ on $\Cl(\cO_XG)$ is a homomorphism.

{\bf Proof}. Since the canonical homomorphism $K_0T(\cO_XG) \ra \Cl(\cO_XG)$
is surjective, it suffices to show the corresponding assertions for $K_0T(\cO_XG)$
(by Lemma 1.4). By Lemma 2.3, we may furthermore assume that $X=\Spec(A)$
where $A$ is a local Dedekind domain. Then, the connecting homomorphism
$\partial: K_1(FG) \ra K_0T(\cO_XG)$ is surjective (see the proof of
Proposition 2.4), and Proposition 2.6 follows from Lemma 1.5.

{\bf Theorem 2.7}. Let $G$ be Abelian and $\gcd(i, \ord(G))=1$. We fix $i' \in
\NN$ such that $ii' \equiv 1 \textrm{ mod } e(G)$ where $e(G)$ denotes the exponent of $G$.
Let $\phi_{i'}$ denote both the
$\cO_X$-algebra automorphism $\cO_XG \ra \cO_XG$ given by $[g] \mapsto [g^{i'}]$ 
and the automorphism
of $K_0(\cO_XG)$ or $\Cl(\cO_XG)$ induced by the association $[\cP] \mapsto
[\cO_XG \otimes_{\cO_XG} \cP]$ (where $\cO_XG$ is considered as an $\cO_XG$-algebra
via $\phi_{i'}$). Then we have:
\[\sigma^i=\phi_{i'} \quad \textrm{on} \quad \Cl(\cO_XG). \]

{\bf Proof}. As in Proposition 2.6, it suffices to show the corresponding 
assertion for $K_1(FG)$ where $\phi_{i'}$ on $K_1(FG)$ is defined analogously.
Since $FG$ is semilocal and commutative, the canonical homomorphism $(FG)^\times \ra
K_1(FG)$ is bijective (see Corollary (9.2) on p.\ 267 in \cite{Ba}). Under this
isomorphism, the automorphism $\phi_{i'}$ corresponds to the restriction of the
(analogously defined) automorphism $\phi_{i'}$ of $FG$. Thus it suffices to
show that the following diagram commutes:
\[\xymatrix{
(FG)^\times \ar[r]^{\sim} \ar[d]^{\phi_{i'}} & K_1(FG) \ar[d]^{\sigma^i} \\
(FG)^\times \ar[r]^{\sim} & K_1(FG). }\]
Now, let $W$ be a local domain of characteristic $0$ whose residue class field
is isomorphic to $F$. (If $\textrm{char}(F) =0$, we may choose $F$ itself for $W$.
If $p =\textrm{char}(F) > 0$, the ring of infinite Witt vectors
over $F$ associated with the prime $p$ is such a ring.) Since the group ring
$WG$ is semilocal and commutative, the canonical map $(WG)^\times \ra K_1(FG)$
is bijective (see loc.\ cit.) and the canonical homomorphism $(WG)^\times \ra
(FG)^\times$ is surjective. Thus it suffices to show that the following diagram 
commutes:
\[\xymatrix{
(WG)^\times \ar[r]^{\sim} \ar[d]^{\phi_{i'}} & K_1(WG) \ar[d]^{\sigma^i}\\
(WG)^\times \ar[r]^\sim & K_1(WG).}\]
In a similar way, we conclude that it suffices to show that the corresponding
diagram commutes if $W$ is replaced by the quotient field $Q$ of $W$ and finally
by the algebraic closure $\bar{Q}$ of $Q$. In the latter case, the 
commutativity follows from Theorem 1.6(d) in \cite{KoCl}, Theorem 3.3 in
\cite{KoAdHi}, and Lemma 3.6(b) in \cite{KoCl}. This ends the proof of
Theorem 2.7.

{\bf Remark 2.8}. Let $\gcd(i,\ord(G))=1$. Theorem 2.7 implies in particular
that $\sigma^{i + e(G)} = \sigma^i$ on $\Cl(\cO_XG)$ if $G$ is Abelian. This
also holds if $X= \Spec(\cO_F)$ where $\cO_F$ is the ring of integers in a
number field $F$ (see Corollary 3.8 in \cite{KoCl}) or if $X$ is a smooth
curve over a finite field (this follows from Theorem 2.10). It is not
clear to me whether this is true in general.

Now, let $L$ be an algebraically closed field such that $\textrm{char}(L)$
does not divide $\ord(G)$, and let $p: X\ra \Spec(L)$ be an irreducible 
smooth curve over $L$. Then, for any f.\ g.\ $LG$-module $V$, 
the pull-back $p^*(V)$ is a locally projective coherent $\cO_XG$-module.
Furthermore, for any locally projective coherent $\cO_XG$-module $\cP$, $\cP'$,
the $\cO_X$-module ${\cal H}{\rm om}_{\cO_X}(\cP,\cP') \cong \cP^\vee \otimes_{\cO_X} \cP$
is again a locally projective $\cO_XG$-module. Finally, for any locally projective
$\cO_XG$-module $\cP$, the $\cO_X$-module $\cP^G$ of $G$-fixed elements is locally
free since $\ord(G)$ is invertible on $X$. Thus, we obtain a well-defined
homomorphism
\[\xymatrix@R=1ex{
K_0(\cO_XG) \ar[r] & \Hom(K_0(LG), K_0(X))\\
[\cP] \ar@{|->}[r]& ([V] \mapsto [\Hom_{\cO_XG}(p^*(V), \cP)]). }\]
This homomorphism is bijective (see the proof of Proposition (2.2) on p.\ 133 in
\cite{Se}) and induces an isomorphism
\begin{equation}
\Cl(\cO_XG) \,\, \tilde{\ra} \,\, \Hom(K_0(LG), \Cl(X))
\end{equation}
by Proposition 2.4.

Let $\psi^i$ denote the $i$-th Adams operation on $K_0(LG)$. In the sequel,
we will identify $K_0(LG)$ with the ring of virtual characters of $G$. Then 
$\psi^i$ maps a character $\chi$ to the character $G \ra L$, $g \mapsto
\chi(g^i)$. Let $\hat{\psi}^i$ denote the adjoint operation (with respect to the
usual character pairing). Note that the assumption $\textrm{char}(L) \notteilt
\ord(G)$ implies that $\gcd(i,\ord(G))$ is invertible on $G$ for all $i\in \NN$.

{\bf Theorem 2.9}. Under the isomorphism (\theequation), the operation $\sigma^i$ on
$\Cl(\cO_XG)$ corresponds to the endomorphism $\Hom(\hat{\psi}^i, \Cl(X))$ of
$\Hom(K_0(LG), \Cl(X))$. 

{\bf Proof}. By Theorem 3.3 on p.\ 145 in \cite{KoAdHi} and Theorem 1.6(d)(ii) in
\cite{KoCl}, 
the operation $\sigma^i$ on $K_1(FG)$ (constructed e.g.\ in section 1) corresponds
to the endomorphism $\Hom(\hat{\psi}^i, K_1(F))$ of $\Hom(K_0(LG), K_1(F))$ 
under the isomorphism 
\[\xymatrix@R=1ex{
K_1(FG) \ar[r]^\sim & \Hom(K_0(LG), K_1(F)) \\
(P,\alpha) \ar@{|->}[r] & 
\Big([V] \mapsto (\Hom_{FG}(F\otimes_L V,P), \Hom_{FG}(F\otimes_L V, \alpha))\Big).}\]
For
any closed point $x\in X$, the association $[M] \mapsto ([V] \mapsto 
[\Hom_{\cO_{X,x}G}(\cO_{X,x} \otimes_L V, M)])$ induces an isomorphism
$K_0T(\cO_{X,x}G) \,\, \tilde{\ra}\,\, \Hom(K_0(LG), K_0T(\cO_{X,x}))$ (both
sides are isomorphic to $K_0(LG)$!) such that the following diagram commutes:
\[\xymatrix{
K_1(FG) \ar@{->>}[rrrr]^\partial \ar[d]^\wr &&&& K_0T(\cO_{X,x}G) \ar[d]^\wr\\
\Hom(K_0(LG), K_1(F)) \ar@{->>}[rrrr]^{\Hom(K_0(LG), \partial)} &&&& 
\Hom(K_0(LG),K_0T(\cO_{X,x})).}\]
Hence, by Lemma 1.5, the operation $\sigma^i$ on $K_0T(\cO_{X,x}G)$ corresponds
to the endomorphism $\Hom(\hat{\psi}^i, K_0T(\cO_{X,x}))$ of 
$\Hom(K_0(LG),K_0T(\cO_{X,x}))$. Under the isomorphism of Lemma 2.3, the
operation $\sigma^i$ on $K_0T(\cO_{X}G)$ obviously corresponds to the 
endomorphism $\oplusm_{x\in X \;{\rm closed}} \sigma^i$ of 
$\oplusm_{x\in X \;{\rm closed}} K_0T(\cO_{X,x}G)$. Thus, under the isomorphism
$K_0T(\cO_XG) \cong \Hom(K_0(LG), K_0T(\cO_X))$, $[\cM] \mapsto 
([V] \mapsto [\Hom_{\cO_XG}(p^*(V),\cM)])$, the operation $\sigma^i$ on
$K_0T(\cO_XG)$ corresponds to the endomorphism $\Hom(\hat{\psi}^i, K_0T(\cO_X))$
of $\Hom(K_0(LG),K_0T(\cO_X))$. Furthermore, the following diagram obviously
commutes:
\[\xymatrix{
K_0T(\cO_XG) \ar[r]^{{\rm can}} \ar[d]^\wr & K_0(\cO_XG) \ar[d]^\wr \\
\Hom(K_0(LG), K_0T(\cO_X)) \ar[r]^{{\rm can}} & \Hom(K_0(LG), K_0(X)).}\]
Now, Theorem 2.9 follows from Lemma 1.4 and Proposition 1.6.

Now, let $L$ be a finite field with $\textrm{char}(L) \notteilt \ord(G)$ and 
$p: X\ra \Spec(L)$ an irreducible smooth curve over $L$. Let 
$\bar{L}$ denote an algebraic closure of $L$ and $\bar{p}: \bar{X}
:= X\times_L \bar{L} \ra \Spec(\bar{L})$ the corresponding curve over $\bar{L}$. 
Then, the composition of the canonical map $K_0(\cO_XG) \ra K_0(\cO_{\bar{X}}G)$
with the isomorphism $K_0(\cO_{\bar{X}}G) \cong \Hom(K_0(\bar{L}G), K_0(\bar{X}))$
constructed above obviously induces a homomorphism
\[K_0(\cO_XG) \ra \Hom_{{\rm Gal}(\bar{L}/L)}(K_0(\bar{L}G), K_0(\bar{X})).\]

{\bf Theorem 2.10}. This homomorphism is bijective. In particular, we obtain
an isomorphism 
\[\Cl(\cO_XG) \,\, \tilde{\ra} \,\, \Hom_{{\rm Gal}(\bar{L}/L)}
(K_0(\bar{L}G), \Cl(\bar{X})).\]
Under this isomorphism, the operation $\sigma^i$ on $\Cl(\cO_XG)$ corresponds
to the endomorphism $\Hom_{{\rm Gal}(\bar{L}/L)}(\hat{\psi}^i, \Cl(\bar{X}))$
of $\Hom_{{\rm Gal}(\bar{L}/L}(K_0(\bar{L}G), \Cl(\bar{X}))$.

{\bf Proof}. The bijectivity can be shown as in section 6 of \cite{ABGr} using
Morita equivalence and the Galois descent property $K_0(X \times_L L') \cong
K_0(\bar{X})^{{\rm Gal}(\bar{L}/L')}$ (for any finite extension $L \subseteq L'
\subset \bar{L}$ of $L$). Proposition 2.4 then yields the Hom-description of
the classgroup. The last assertion immediately follows from Theorem 2.9.

\bigskip

\section*{\S 3 Equivariant Riemann-Roch Type Formulas for Tame Extensions of
Dedekind Schemes}

The aim of this section is to prove Theorem A and Theorem B presented in the introduction.

Let $Y$ be a Dedekind scheme and $G$ a finite group of order $n$. Let 
$\Ind_1^G: \Cl(\cO_Y) \ra \Cl(\cO_YG)$ and $\Ind_1^G: K_0T(\cO_Y) \ra
K_0^\lf T(\cO_YG)$ denote the induction maps. The following lemma generalizes
Lemma 2.6 on p.\ 933 in \cite{BC}. 

{\bf Lemma 3.1}. The image of the natural multiplication maps 
\[ K_0T(\cO_Y) \times K_0^\lf(\cO_YG) \ra K_0^\lf T(\cO_YG) \quad 
\textrm{and} \quad \Cl(\cO_Y) \times K_0^\lf(\cO_YG) \ra \Cl(\cO_YG)\]
is contained in $\Ind_1^GK_0T(\cO_Y)$ resp.\ $\Ind_1^G\Cl(\cO_Y)$.

{\bf Proof}. The assertion for the first
map is clear. The assertion for the second map follows from this since the
natural map $K_0T(\cO_Y) \ra \Cl(\cO_Y)$ is surjective.

Now, let $F/E$ be a finite Galois extension of the function field $E$ of $Y$
with Galois group $G$. Let $X$ denote the normalization of $Y$ in $F$. Then
$X$ is a Dedekind scheme endowed with a natural $G$-action and the 
corresponding $G$-morphism $f:X\ra Y$ is finite (see the proof of Theorem (8.1)
on p.\ 47 in \cite{N}). We assume that $f$ is tamely ramified. As in Lemma 5.5
in \cite{KoCl}, one easily shows that then, for any locally free coherent
$\cO_X$-module $\cE$ with (semilinear) $G$-action, the direct image $f_*(\cE)$
is a locally free coherent $\cO_YG$-module in the sense of section 2. Let
$K_0(G,X)$ denote the Grothendieck group of all such modules $\cE$. Thus,
we have a homomorphism
\[f_*: K_0(G,X) \ra K_0^\lf(\cO_YG), \quad [\cE] \mapsto [f_*(\cE)].\]
The different $\cD:= \cD_{X/Y} := {\cal A}{\rm nn}_{\cO_X}(\Omega^1_{X/Y})$
is a $G$-stable ideal in $\cO_X$, hence a module $\cE$ as above. The
following proposition generalizes formula (2.8) on p.\ 933 in \cite{BC}.

{\bf Proposition 3.2}. For all $x\in K_0(G,X)$ we have:
\[f_*\left(x \cdot \sum_{i=0}^{n-1}[\cD^{-i}]\right) =0 \quad \textrm{in} \quad
K_0^\lf(\cO_YG)/(\Ind_1^G\Cl(\cO_Y) \oplus n \ZZ [\cO_YG]).\]

{\bf Proof}. We may assume that $x=[\cE]$ where $\cE$ is a module as above.
Let $r:= \rank_{\cO_X}(\cE)$. Then we have:
\begin{eqnarray*}
\lefteqn{\sum_{i=0}^{n-1}\Big([f_*(\cE\otimes \cD^{-i})] - r [\cO_YG]\Big) }\\
&=&n\Big([f_*(\cE)] - r[\cO_YG]\Big) + \sum_{i=1}^{n-1}\Big([f_*(\cE \otimes \cD^{-i})]
- [f_*(\cE)]\Big) \quad \textrm{in} \quad \Cl(\cO_YG).
\end{eqnarray*}
In the sequel, let $\cM \mapsto \cM^t$ denote the forgetful functor from the
category of $\cO_YG$-modules to the category of $\cO_Y$-modules. (We will
consider $\cM^t$ also as an $\cO_YG$-module with trivial $G$-action.) Then, the
elements $[f_*(\cO_X)^t]-n$ and $[f_*(\cE)^t] -nr$ are contained in $\Cl(\cO_Y)$.
Hence, we have by Lemma 3.1:
\[n([f_*(\cE)] -r [\cO_YG]) = [f_*(\cO_X)^t \otimes f_*(\cE)] -
[f_*(\cE)^t \otimes \cO_YG] \quad \textrm{in} \quad 
\Cl(\cO_YG)/\Ind_1^G\Cl(\cO_Y).\]
The homomorphism 
\[f_*(\cO_X)^t \otimes f_*(\cE) \ra f_*(\cE)^t \otimes \cO_YG, 
\quad a\otimes b \mapsto \sum_{g\in G} a g(b) \otimes [g^{-1}],\]
of $\cO_YG$-modules is generically bijective since $F/E$ is a Galois extension
and any f.\ g.\ module over the twisted group ring $F\#G$ is isomorphic to
$\oplusm^m F$ for some $m\ge 0$. In particular, this map is a monomorphism
and the cokernel $\cR_{X/Y}(\cE)$ is an $\cO_YG$-torsion module. Hence, it suffices
to show that we have:
\[[\cR_{X/Y}(\cE)] = \sum_{i=1}^{n-1}[f_*(\cE\otimes \cD^{-i}/\cO_X)] \quad
\textrm{in} \quad K_0^\lf T(\cO_YG)/\Ind_1^GK_0T(\cO_Y).\]
By Lemma 2.3, it furthermore suffices to show that we have 
\[[\cR_{X/Y}(\cE)_y] = \sum_{i=1}^{n-1}[f_*(\cE \otimes \cD^{-i}/\cO_X)_y] \quad
\textrm{in} \quad K_0^\lf T(\cO_{Y,y}G)/\Ind_1^GK_0T(\cO_{Y,y})\]
for all closed points $y \in Y$. \\
We now fix $y\in Y$ and $x\in X$ with $f(x) =y$. Let $G_x := \{g\in G: xg =x\}$
denote the decomposition group of $x$. Furthermore, let $f':X':= \Spec(\hat{\cO}_{X,x})
\ra \Spec(\hat{\cO}_{Y,y}) =: Y'$ denote the induced $G_x$-morphism where $\hat{}$ 
denotes completion. We identify the category of coherent torsion modules on $Y'$
with the category of coherent torsion modules on $Y$ supported in $y$. 
An easy generalization of Corollary 3.11(b) on p.\ 239 in
\cite{Ch} shows that $\cR_{X/Y}(\cE)_y$ is isomorphic to the direct sum of 
$[G:G_x]$ copies of $\Ind_{G_x}^G \cR_{X'/Y'}(\hat{\cE}_x)$. Furthermore, it is
clear that $f_*(\cE \otimes \cD^{-i}_{X/Y}/\cO_X)_y$ is isomorphic to
$\Ind_{G_x}^G f'_*(\hat{\cE}_x \otimes \cD^{-i}_{X'/Y'}/\cO_{X'})$ for all 
$i \ge 0$. For $i \equiv j \textrm{ mod } \ord(G_x)$, we finally have
\[[f'_*(\hat{\cE}_x \otimes \cD_{X'/Y'}^{-i}/\cO_{X'})] =
[f'_*(\hat{\cE}_x \otimes \cD^{-j}_{X'/Y'}/\cO_{X'})] \quad \textrm{in}  \quad
K_0^\lf T(\cO_{Y'}G_{x})/\Ind_1^{G_x} K_0T(\cO_{Y'})\]
since the ideal $\cD^{\ord(G_x)}_{X'/Y'}$ of $\cO_{X'}$ can be written 
as $(f')^*({\mathfrak a})$ with some ideal ${\mathfrak a}$ in $\cO_{Y'}$ and 
since, for any locally free coherent $\cO_{Y'}G$-module $\cP$, we have
\[[\cP/{\mathfrak a} \cP] = [\cO/{\mathfrak a} \otimes \cP] = 0 \quad
\textrm{in} \quad K_0^\lf(\cO_{Y'}G_{x})/\Ind_1^{G_x} K_0T(\cO_{Y'})\]
by Lemma 3.1. Thus it suffices to prove that
\[[\cR_{X'/Y'}(\hat{\cE}_x)] = \sum_{i=1}^{\ord(G_x)-1} 
[f'_*(\hat{\cE}_x \otimes \cD^{-i}_{X'/Y'} / \cO_{X'})] \quad \textrm{in}
\quad K_0T(\cO_{Y'}G_x)/\Ind_1^{G_x}K_0T(\cO_{Y'}).\]
We now write $G$ for $G_x$, $X$ for $X'$, $\cE$ for $\hat{\cE}_x$, and so on. Let 
$\Delta \subseteq G$ denote the inertia group, $e$ the order of $\Delta$,
${\mathfrak P}$ the ideal in $\cO_X$ which corresponds to the closed point in $X$,
and $\chi$ the $\Delta$-module ${\mathfrak P}/{\mathfrak P}^2$. We decompose
$f: X\ra Y$ into $X\,\, \stackrel{g}{\ra}\,\, Z \,\, \stackrel{h}{\ra} \,\,Y$
where $Z:= \Spec(\Gamma(X,\cO_X)^\Delta)$; i.e., the function field of $Z$ is
the inertia field of $F/E$. Since $K_0(G,X)$ is generated by the classes
of fractional $G$-stable ideals in $\cO_X$ (see Lemma 5.5(c) in \cite{KoCl}),
we may assume that $\cE = {\mathfrak P}^j$ for some $j\in \ZZ$. An easy
generalization of Corollary 3.8 on p.\ 236 and Theorem 2.8 on p.\ 222 in \cite{Ch}
shows that we have the following isomorphisms:
\begin{eqnarray*}
\lefteqn{\cR_{X/Y}(\wP^j) \cong \Ind_\Delta^G h_*(\cR_{X/Z}(\wP^j))}\\
&\cong&\Ind_\Delta^Gh_*\left(\oplusm_{i=1}^{e-1} g_*\Big((\wP^j/\wP^{j+i})^t \otimes
\chi^{j+i}\Big)\right) \\
&\cong& \Ind_\Delta^G f_*\left(\oplusm_{i=1}^{e-1}(\wP^j/\wP^{j+i})^t \otimes  
\chi^{j+i} \right).
\end{eqnarray*}
Thus we have:
\[[\cR_{X/Y}(\wP^j)] = \sum_{i=1}^{e-1} i [\Ind_\Delta^G f_*(\chi^{j+i})] \quad
\textrm{in} \quad K_0T(\cO_YG).\]
Since $\cD= \wP^{e-1}$ and $\wP^e = f^*(\wp)$ (where $\wp$ is the ideal in
$\cO_Y$ which corresponds to the closed point in $Y$), we can conclude as
above using Lemma 3.1:
\begin{eqnarray*}
\lefteqn{\sum_{i=1}^{n-1}[f_*(\wP^j\otimes \cD^{-i}/\cO_X)] =
\frac{n}{e} \sum_{i=1}^{e-1}[f_*(\wP^j \otimes \cD^{-i}/\cO_X)] =
\frac{n}{e} \sum_{i=1}^{e-1} [f_*(\wP^{j+i}/\wP^{j+e})]}\\
&=& \frac{n}{e} \sum_{i=1}^{e-1} i [f_*(\wP^{j+i}/\wP^{j+i+1})] \quad
\textrm{in} \quad K_0T(\cO_YG)/\Ind_1^GK_0T(\cO_Y). \hspace*{10ex}
\end{eqnarray*}
Thus it suffices to prove that the $\cO_YG$-modules $\Ind_\Delta^G f_*(\chi^i)$
and $\oplusm^{n/e}f_*(\wP^i/\wP^{i+1})$ are isomorphic for all $i\in \ZZ$. For this,
we consider the $\cO_YG$-homomorphism
\[\xymatrix@R=1ex{
h_*(\cO_Z)^t \otimes f_*(\wP^i/\wP^{i+1})
\ar[r]& {\rm Maps}_\Delta(G, f_*(\wP^i/\wP^{i+1})) \\
a\otimes b \ar@{|->}[r] & (g \mapsto ag(b)).}\]
This homomorphism is bijective since $h$ is unramified (e.g., see pp.\ 214-215 in
\cite{Ch}). Furthermore, the left hand side is obviously isomorphic to
$\oplusm^{n/e}f_*(\wP^i/\wP^{i+1})$ and the right hand side is isomorphic to
$\Ind_\Delta^G f_*(\chi^i)$. So, Proposition 3.2 is proved.

Now, let $k\in \NN$ with $\gcd(k,n)=1$ and $k'\in \NN$ with $kk' \equiv 1
\textrm{ mod } n$. Let $\sigma^k$ denote the $k$-th symmetric power operation
on $K_0(G,Y)$ and $\psi^k$ the $k$-th Adams operation on $K_0(G,Y)$ or 
$K_0(G,X)$ (e.g., see section 1 in \cite{KoCl}). The composition of the map $f_*:
K_0(G,X) \ra K_0^\lf (\cO_YG)$ with the Cartan homomorphism $K_0^\lf(\cO_YG)
\ra K_0(G,Y)$ is denoted by $f_*$ again. Finally, let $\hat{K}_0(G,Y)[k^{-1}]$
denote the $J$-adic completion of $K_0(G,Y)[k^{-1}]$ where 
$J:={\rm ker}(K_0(G,Y) \,\, \stackrel{{\rm rank}}{\longrightarrow}\,\, \ZZ)[k^{-1}]$
is the augmentation ideal in $K_0(G,Y)[k^{-1}]$. 

{\bf Theorem 3.3}. For all $x \in K_0(G,X)$ we have:
\[\sigma^k(f_*(x)- \rank(x)\cdot [\cO_YG]) = 
f_*\left(\sum_{i=0}^{k'-1}[\cD^{-ik}] \cdot \psi^k(x)\right)\]
in $\hat{K}_0(G,Y)[k^{-1}]/(\Ind^G_1 K_0(Y)) \hat{K}_0(G,Y)[k^{-1}]$.

{\bf Proof}. Let 
\[\hat{f}_*: \hat{K}_0(G,X)[k^{-1}] := K_0(G,X) \otimes_{K_0(G,Y)} 
\hat{K}_0(G,Y)[k^{-1}] \ra \hat{K}_0(G,Y)[k^{-1}]\]
denote the homomorphism which is induced by $f_*: K_0(G,X) \ra K_0(G,Y)$, and
let $\theta^k(\cD^{-1}):= 1+[\cD^{-1}] + \ldots + [\cD^{-(k-1)}] \in K_0(G,X)$ 
denote the Bott element. As in Theorem 5.4 in \cite{KoCl}, one easily deduces
the following assertion from the equivariant Adams-Riemann-Roch theorem (see
Theorem (4.5) in \cite{KoGRR}): The element $\theta^k(\cD^{-1})$ is invertible
in $\hat{K}_0(G,X)[k^{-1}]$ and we have
\[\psi^k(f_*(x)) = \hat{f}_*(k \cdot \theta^k(\cD^{-1})^{-1} \cdot \psi^k(x))
\quad {\rm in} \quad \hat{K}_0(G,Y)[k^{-1}]\]
for all $x\in K_0(G,X)$.
Furthermore, we have:
\[\theta^k(\cD^{-1}) \cdot \left(\sum_{i=0}^{k'-1}[\cD^{-ik}]\right) =
\sum_{j=0}^{k-1}\sum_{i=0}^{k'-1}[\cD^{-(j+ik)}] =
\sum_{i=0}^{kk'-1} [\cD^{-i}] = [\cO_X] + \sum_{i=1}^{kk'-1}[\cD^{-i}]\]
in $K_0(G,X)$. Thus, we have:
\[\theta^k(\cD^{-1})^{-1} = \sum_{i=0}^{k'-1}[\cD^{-ik}] -
\theta^k(\cD^{-1})^{-1}\sum_{i=1}^{kk'-1} [\cD^{-i}] \quad 
\textrm{in} \quad \hat{K}_0(G,X)[k^{-1}].\]
Hence, we obtain the equality
\begin{eqnarray*}
\lefteqn{\psi^k(f_*(x)) = k \cdot \hat{f}_*\left(\left( \sum_{i=0}^{k'-1}[\cD^{-ik}] -
\theta^k(\cD^{-1})^{-1}\cdot \sum_{i=1}^{kk'-1}[\cD^{-i}]\right) \cdot
\psi^k(x)\right)}\\
&=& k \cdot f_*\left(\sum_{i=0}^{k'-1}[\cD^{-ik}]\cdot \psi^k(x)\right) \quad
\textrm{in} \quad \hat{K}_0(G,Y)[k^{-1}]/(\Ind_1^GK_0(Y)) \hat{K}_0(G,Y)[k^{-1}]
\end{eqnarray*}
by Proposition 3.2. Since we have $\psi^k = k \cdot \sigma^k$ on $\Cl(\cO_YG)$ (by
Proposition 2.6) and $\psi^k([\cO_YG]) = [\cO_YG]$ (by Theorem 1.6(e) in
\cite{KoCl}), this implies Theorem 3.3.

Note that the formula of Theorem 3.3 lives within the somewhat complicated
group $\hat{K}_0(G,Y)[k^{-1}]/(\Ind_1^GK_0(Y)) \hat{K}_0(G,Y)[k^{-1}]$.
The next proposition computes this group in a special case.

{\bf Proposition 3.4}. Let $L$ be an algebraically closed field, $Y$ a projective
smooth irreducible curve over $L$, and $n= \ord(G)$ a power of a prime $l\not=
{\rm char}(L)$. Let $I$ denote the augmentation ideal in $K_0(LG)$. Then we have:
\[\hat{K}_0(G,Y)[k^{-1}] \cong K_0(Y)[k^{-1}] \oplus I\otimes \ZZ_l \oplus
I\otimes \ZZ_l;\]
under this isomorphism, the extended ideal $(\Ind_1^G K_0(Y))\hat{K}_0(G,Y)[k^{-1}]$ 
corresponds to the subgroup $\{(ny,([\ZZ G]-n) \otimes \rank(y), ([\ZZ G] -n) \otimes
\deg \det(y)): y \in K_0(Y)[k^{-1}]\}$ of $K_0(Y)[k^{-1}] \oplus 
I \otimes \ZZ_l \oplus I \otimes \ZZ_l$. 

{\bf Proof}. The canonical map $K_0(LG) \otimes K_0(Y) \ra K_0(G,Y)$ is an
isomorphism by Proposition (2.2) on p.\ 133 in \cite{Se}. Since the augmentation ideal
of $K_0(Y)$ is nilpotent (e.g., by Proposition 2.6) and the $I$-adic topology
on $I$ coincides with the $l$-adic topology (see Proposition 1.1 on p.\ 277
in \cite{AT}), the completion $\hat{K}_0(G,Y)[k^{-1}]$ is isomorphic to the 
direct sum of $K_0(Y)[k^{-1}]$ and the $l$-adic completion of $I\otimes K_0(Y)[k^{-1}]$.
Furthermore, we have $K_0(Y) \cong \ZZ \oplus \ZZ \oplus \Pic^0(Y)$ where 
$\Pic^0(Y)$ denotes the group of line bundles on $Y$ of degree $0$. Since 
$\Pic^0(Y)$ is an $l$-divisible group (see item (iv) on p.\ 42 in \cite{Mu}),
the $l$-adic completion of $I\otimes K_0(Y)[k^{-1}]$ is isomorphic to
$I\otimes \ZZ_l \oplus I\otimes \ZZ_l$. Thus, we have
\[\hat{K}_0(G,Y)[k^{-1}] \cong K_0(Y)[k^{-1}] \oplus I\otimes \ZZ_l \oplus
I\otimes \ZZ_l.\]
Under the isomorphism $K_0(G,Y) \cong K_0(LG) \otimes K_0(Y)$, the ideal
$\Ind_1^G K_0(Y)$ of the ring $K_0(G,Y)$ 
corresponds to the ideal $\Ind_1^G K_0(L) \otimes
K_0(Y)$ ($\cong K_0(Y)$) of $K_0(LG) \otimes K_0(Y)$ which is generated by the
element $[\ZZ G] \otimes 1 = n\otimes 1 + ([\ZZ G]-n) \otimes 1$. One easily
deduces the second assertion of Proposition 3.4 from this. (Note that
$[\ZZ G] \cdot x =0$ for all $x \in I$.)

Now, let $f_*: K_0(G,X) \ra \Cl(\cO_YG)$ denote the composition of 
$f_*:K_0(G,X) \ra K_0^\lf(\cO_YG)$ with the canonical projection
$K_0^\lf(\cO_YG) \cong \Cl(\cO_Y G) \oplus \ZZ[\cO_Y G] \ra \Cl(\cO_Y G)$. 

{\bf Theorem 3.5}. Suppose that one of the following conditions holds:\\
(a) $Y= \Spec(\cO_E)$ where $\cO_E$ is the ring of integers in a number field $E$.\\
(b) $Y$ is an irreducible projective smooth curve over a finite field $L$ and \\
$\gcd({\rm char}(L), n) =1$. \\
(c) The group $G$ is Abelian and $f:X\ra Y$ is unramified.\\
(d) $k=1$.\\
Then we have for all $x\in K_0(G,X)$:
\[\sigma^k(f_*(x)) = f_*\left(\sum_{i=0}^{k'-1} [\cD^{-ik}] \cdot \psi^k(x)\right)
\quad \textrm{in} \quad \Cl(\cO_YG)/\Ind_1^G\Cl(\cO_Y).\]

{\bf Proof}. In the case (a), Theorem 3.5 can be deduced from Corollary 2.7 on p.\ 933 in
\cite{BC} using Theorem 3.7 and Lemma 5.5 in \cite{KoCl} (see also the proof
of Theorem 5.6 in \cite{KoCl}). The same can be done in the case (b) by using 
Lemma 3.6(a) in \cite{KoCl} and
Theorem 2.9 (in place of Theorem 3.7 in \cite{KoCl}) and an obvious generalization
of Lemma 5.5(c) in \cite{KoCl}. (For completeness sake, we mention that it
is easy to check that the additional assumptions in Theorem 2.1 on p.\ 932 
in \cite{BC} about the absolute discriminant or the characteristic of $E$ are 
not necessary for Corollary 2.7 on p.\ 933 in \cite{BC}.) We now prove Theorem 3.5
in the case~(c), i.e., we want to show the formula
\begin{equation}
\sigma^k(f_*(x)) = k'\cdot f_*(\psi^k(x)) \quad \textrm{in} \quad 
\Cl(\cO_YG)/\Ind_1^G\Cl(\cO_Y)
\end{equation}
for all $x\in K_0(G,X)$. First, we show that it suffices to prove the formula
(\theequation) for $x=1=[\cO_X]$.  Indeed, for an arbitrary $x \in K_0(G,X)$,
there is a $y \in K_0(Y) \subseteq K_0(G,Y)$ such that $x = f^*(y)$ (e.g., see
Theorem 1(B) on p.\ 112 in \cite{Mu}). Furthermore, we have 
$\sigma^k(\Ind_1^G\Cl(\cO_Y)) \subseteq \Ind_1^G\Cl(\cO_Y)$. This follows
from Proposition 1.1 in \cite{KoCl} as there is a polynomial $Q_k \in
\ZZ[X_1, \ldots, X_k; Y_1, \ldots, Y_k]$ which is homogeneous of weight $k$
in both sets of variables such that 
\[\sigma^k(z \cdot [\cO_Y G]) = Q_k\left(\sigma^1(z), \ldots, \sigma^k(z); 
[\Sym^1(\cO_YG)], \ldots, [\Sym^k(\cO_YG)]\right)\;\; \textrm{in} \;\;
\Cl(\cO_YG)\]
for all $z \in Cl(\cO_Y)$ (by Theorem 2.2 in \cite{KoCl}). Thus we have
in $\Cl(\cO_YG)/\Ind_1^G\Cl(\cO_Y)$:
\begin{eqnarray*}
\lefteqn{\sigma^k(f_*(x)) = \sigma^k(f_*(f^*(y))) = \sigma^k(y \cdot f_*(1))
\quad \textrm{(Projection formula)}}\\
&=& \sigma^k(\rank(y) \cdot f_*(1)) \quad \textrm{(Lemma 3.1)}\\
&=& \rank(y) \cdot \sigma^k(f_*(1)) \quad \textrm{(Proposition 2.6)}\\
&=& \rank(y) \cdot k' \cdot f_*(1) \quad \textrm{(by assumption)}\\
&=& k' \cdot \psi^k(y) \cdot f_*(1) \quad \textrm{(Lemma 3.1)}\\
&=& k' \cdot f_*(\psi^k(f^*(y))) = k' \cdot f_*(\psi^k(x)) \quad
\textrm{(Projection formula).}
\end{eqnarray*}
We now prove formula (\theequation) for $x=1$. Since $f$ is unramified, 
the scheme $X$ is a principal $G$-bundle over $Y$ (see Proposition 2.6 on p.\
115 in \cite{SGA1}). There is a well-known natural bijection between the
set of all principal $G$-bundles over $Y$ and the cohomology group $H^1(Y,G)$.
We write $[X]$ for the corresponding element in $H^1(Y,G)$. We define
a new principal $G$-bundle $X_{k'}$ over $Y$ as follows: $X_{k'} =X$ as
$Y$-schemes and the new action $*$ of $G$ on $X_{k'}$ is given by
$x*g:=xg^k$ for ``$x \in X$'' and $g\in G$. Then, it is easy to check that
the association $X \mapsto X_{k'}$ corresponds to the 
multiplication with $k'$ on $H^1(Y,G)$. 
Let ${\rm cl}: H^1(Y,G) \ra \Cl(\cO_Y G)$ denote the map which maps a
principal $G$-bundle $f: X \ra Y$ to the class $[f_*(\cO_X)] - [\cO_YG]$. 
This map is a homomorphism by Theorem 5 and the subsequent remarks on p.\ 189
in \cite{Wa} and by Proposition 3.9 in \cite{ABGr}. Thus
we have:
\begin{eqnarray*}
\lefteqn{\sigma^k(f_*([\cO_X])) = \phi_{k'}({\rm cl}([X])) \quad
\textrm{(Theorem 2.7)}}\\
&=& {\rm cl}([X_{k'}]) = {\rm cl}(k' \cdot [X]) 
= k' \cdot {\rm cl}([X]) = k' \cdot f_*([\cO_X]).
\end{eqnarray*}
in $\Cl(\cO_Y G)$, as was to be shown. In the case (d), Theorem 3.5 immediately
follows from Proposition 3.2.

{\bf Remark 3.6}. If one of the conditions (a), (b), (c), (d) of Theorem 3.5 is
satisfied, then Theorem 3.3 follows from Theorem 3.5 by passing from 
$\Cl(\cO_YG) \subset K_0(\cO_YG)$ to $K_0(G,Y)$ via Cartan homomorphism and
finally by passing from $K_0(G,Y)$ to the completion $\hat{K}_0(G,Y)[k^{-1}]$ of
$K_0(G,Y)[k^{-1}]$. In particular, in the case (a), the formula of Theorem 3.3
is substantially weaker than the formula in Theorem 3.5, as already the
passage from $K_0(\cO_YG)$ to $K_0(G,Y)$ loses much information. On the other
hand, in the case (b), the formula of Theorem 3.5 modulo torsion follows
from the formula in Theorem 3.3 if $n$ is a power of a prime. This can
be proved as follows. The Cartan homomorphism $K_0(\cO_YG) \ra K_0(G,Y)$ is
bijective since $n$ is invertible on $Y$. Furthermore, the canonical map
$\Cl(\cO_YG)/\Ind_1^G\Cl(\cO_Y) \subseteq K_0(\cO_YG)/\Ind_1^G K_0(Y) \ra
K_0(\cO_{\bar{Y}}G)/\Ind_1^GK_0(\bar{Y})$ is injective by Theorem 2.10. (Here,
$\bar{Y}$ denotes the curve $Y \times_L \bar{L}$ over the algebraic closure
$\bar{L}$ of $L$.) Hence, it suffices to prove the formula
\begin{equation}
\sigma^k(\bar{f}_*(x)) = \bar{f}_*\left(\sum_{i=0}^{k'-1}[\cD_{\bar{X}/\bar{Y}}^{-ik}]
\cdot \psi^k(x)\right) \quad \textrm{in}\quad K_0(G,\bar{Y})_\QQ/
\Ind_1^GK_0(\bar{Y})_\QQ
\end{equation}
for all $x\in K_0(G,\bar{X})$. Furthermore, we have $K_0(G,\bar{Y})_\QQ \cong
K_0(\bar{L}G)_\QQ \otimes K_0(\bar{Y})_\QQ$ and $K_0(G,\bar{Y})_\QQ/\Ind_1^GK_0(\bar{Y})_\QQ
\cong I\otimes K_0(\bar{Y})_\QQ \cong I_\QQ \oplus I_\QQ$ (see the proof of 
Proposition 3.4). On the other hand, $\left(\hat{K}_0(G,\bar{Y})[k^{-1}]/
(\Ind_1^GK_0(\bar{Y}))\hat{K}_0(G,\bar{Y})[k^{-1}]\right)_\QQ$ is isomorphic to
$I\otimes \QQ_l \oplus I\otimes \QQ_l$ by Proposition 3.4. Hence, the canonical
map $K_0(G,\bar{Y})_\QQ/\Ind_1^GK_0(\bar{Y})_\QQ$ $\ra \left(\hat{K}_0(G,\bar{Y})[k^{-1}]/
(\Ind_1^GK_0(\bar{Y}))\hat{K}_0(G,\bar{Y})[k^{-1}]\right)_\QQ$ is injective, and
formula (\theequation) follows from Theorem 3.3.

{\bf Remark 3.7}. \\
(a) If one of the conditions (a), (b), or (c) holds, 
Theorem 3.5 can be slightly strengthened:
It suffices to assume that $k'$ is an inverse modulo the exponent of $G$ (see
\cite{BC} and Theorem 2.7, respectively). It is not clear to me whether this
is true also in the case~(d).\\
(b) Let $Y$ be an irreducible smooth projective curve over a finite field $L$.
Then, the case~(c) is particularly interesting as complementary case of the
semisimple case which is assumed in the case~(b). Indeed, if $G$ is an (Abelian)
${\rm char}(L)$-group, then the tameness condition already implies that $f$
is unramified.

\bigskip

Mathematisches Institut II der Universit\"at Karlsruhe, 
D-76128 Karlsruhe, Germany; {\em e-mail:} 
Bernhard.Koeck@math.uni-karlsruhe.de.


\bigskip


\begin{thebibliography}{SGA\,1}

\bibitem[AB]{ABGr} {\sc A.\ Agboola} and {\sc D.\ Burns}, Grothendieck
groups of bundles on varieties over finite fields, {\em preprint}
(1998), 40 pp.
\bibitem[AT]{AT} {\sc M.\ F.\ Atiyah} and {\sc D.\ O.\ Tall}, 
Group representations, $\lambda $-rings and the $J$-homomorphism, 
{\em Topology} {\bf 8} (1969), 253-297. 
\bibitem[B]{Ba} {\sc H.\ Bass}, ``Algebraic $K$-theory'', {\em Math.\  
Lecture Note Series} (Benjamin, New York, 1968).
\bibitem[BC]{BC} {\sc  D.\ Burns} and {\sc  T.\ Chinburg}, Adams operations 
and integral Hermitian-Galois representations, {\em Amer.\ J.\ Math.}\  
{\bf 118} (1996), 925-962.
\bibitem[CT]{CNT} {\sc  Ph.\ Cassou-Nogu\`es} and {\sc  M.\ J.\ Taylor}, 
Op\'erations d'Adams et Groupe des classes d'Alg\`ebre de groupe, 
{\em J.\ Algebra} {\bf 95} (1985), 125-152.
\bibitem[C]{Ch} {\sc S.\ U.\ Chase}, Ramification invariants and torsion
Galois module structure in number fields, {\em J.\ Algebra} {\bf 91} (1984),
207-257.
\bibitem[CR]{CR} {\sc  C.\ W.\ Curtis} and {\sc  I.\ Reiner}, ``Methods of 
representation theory with applications to finite groups and orders'',
vol.\ I, {\em Pure Appl.\ Math.}\ (Wiley, New York, 1981).
\bibitem[F\,1]{F1} {\sc  A.\ Fr\"ohlich}, ``Galois module structure of algebraic 
integers'', {\em Ergeb.\ Math.\ Grenz\-geb.\ (3)}
{\bf 1} (Springer, New York, 1983).
\bibitem[F\,2]{F2} {\sc  A.\ Fr\"ohlich}, ``Classgroups and Hermitian modules'',
{\em Progr.\ Math.}\  {\bf 48} (Birk\-h\"au\-ser, Boston, 1984).
\bibitem[FL]{FL} {\sc W.\ Fulton} and {\sc S.\ Lang}, 
``Riemann-Roch algebra'', 
{\em Grundlehren Math.\ Wiss.}\ {\bf 277} (Springer, New York, 1985).
\bibitem[G\,1]{Gr1} {\sc  D.\ R.\ Grayson}, Higher algebraic $K$-theory: II,
in {\sc  M.\ R.\ Stein} (ed.), ``Algebraic $K$-theory (Evanston, 1976)'',
{\em Lecture Notes in Math.}\ {\bf 551} (Springer, New York, 1976), 217-240.
\bibitem[G\,2]{Gr2} {\sc  D.\ R.\ Grayson}, Exterior power operations on higher
$K$-theory, {\em $K$-Theory} {\bf 3} (1989), 247-260.
\bibitem[SGA\,1]{SGA1} {\sc A.\ Grothendieck} and {\sc M.\ Raynaud},
``Rev\^etements \'etales et groupe fondamental'', {\em Lecture Notes in
Math.}\ {\bf 224} (Springer, New York, 1971).
\bibitem[K\,1]{KoAdHi} {\sc  B.\ K\"ock}, On Adams operations 
on the higher $K$-theory of group rings, in:
G.\ Banaszak et al.\ (eds.), ``Algebraic $K$-theory (Pozna\'n, 1995)'', 
{\em Contemp.\ Math.}\ {\bf 199} (Amer.\ Math.\ Soc., Providence, 1996), 
139-150.
\bibitem[K\,2]{KoGRR} {\sc B.\ K\"ock}, The Grothendieck-Riemann-Roch
theorem for group scheme actions, {\em Ann.\ Sci.\ 
\'Ecole Norm.\ Sup.} {\bf 31} (1998), 415-458.
\bibitem[K\,3]{KoCl} {\sc B.\ K\"ock}, Operations on locally free 
classgroups, to appear in {\em Math.\ Ann.}\ {\bf 314} (1999), 667-702.
\bibitem[M]{Mu} {\sc D.\ Mumford}, ``Abelian varieties'', {\em Tata Inst.\
of Fundamental Research Stud.\ in Math.} {\bf 5} (Oxford University Press, London,
1970). 
\bibitem[N]{N} {\sc  J.\ Neukirch}, ``Algebraische Zahlentheorie'' (Springer,
Berlin, 1992).
\bibitem[S]{Se} {\sc G.\ Segal}, Equivariant $K$-theory, {\em Publ.\ Math.\
IHES} {\bf 34} (1968), 129-151.
\bibitem[W]{Wa} {\sc W.\ C.\ Waterhouse}, Principal homogeneous spaces and
group scheme extensions, {\em Trans.\ Amer.\ Math.\ Soc.}\ {\bf 153}
(1971), 181-189.


\end{thebibliography}
\end{document}